\documentclass{article}

\input amssym.def
\input amssym.tex

\newcommand{\R}{{\Bbb R}}
\newcommand{\Q}{{\Bbb Q}}
\newcommand{\Z}{{\Bbb Z}}
\newcommand{\N}{{\Bbb N}}
\newcommand{\X}[2]{\mbox{$\hat\chi^{(#1)}(#2)$}}
\newcommand{\Xone}[1]{\mbox{$\hat\chi(#1)$}}
\newcommand{\D}{\Delta}
\renewcommand{\r}{\rho}
\newcommand{\qed}{$\natural\\ \\$}
\newcommand{\colors}{\mbox{$C_1,\ldots,C_m$}}
\newcommand{\by}{\times}

\begin{document}
\addtocounter{footnote}{1}
\title{The $k^{\rm th}$ Upper Chromatic Number of the Line\footnote{This note appeared in \emph{Discrete Mathematics,} vol.~169 (1997), pp.~157--162.}}
\author{Aaron Abrams\cr
University of California, Davis\thanks{Please address all 
correspondence to the author, c/o Mathematics Department, 
University of California, Berkeley CA 94720.}}
\date{}
\maketitle

\newtheorem{thm}{Theorem}
\newtheorem{prop}[thm]{Proposition}
\newtheorem{lem}[thm]{Lemma}
\newtheorem{cor}[thm]{Corollary}
\newtheorem{prob}{Problem}

\begin{abstract}
Let $S \subseteq \R^n$, and let $k\in\N$.  Greenwell and Johnson
\cite{g&j} define \X k S to be the smallest integer $m$ (if such an
integer exists) such that for every $k\by m$ array $D=(d_{ij})$ of
positive real numbers, $S$ can be colored with the colors $\colors$
such that no two points of $S$ which are a (Euclidean) distance $d_{ij}$
apart are
both colored $C_j$, for all $1\leq i \leq k$ and $1\leq j \leq m$.  If no
such integer exists then we say that $\X k S =\infty$.  In this paper we
show that $\X k \R$ is finite for all $k$.
\end{abstract}

\section{Introduction}


Let $S \subseteq \R^n, k,m \in \N$, and let $D=(d_{ij})$ be a $k \by m$
array of positive real numbers.  If there is a coloring $T$ of $S$ using
the colors $\colors$ such that no two points of $S$ that are a distance
$d_{ij}$ apart are both colored $C_j$ (for all $1\leq i \leq k$ and $1\leq j
\leq m$), then we say that $S$ is {\em $D$-colorable}, and we call
$T:S \to \{\colors\}$ a {\em $D$-coloring} of $S$.  We call each
$d_{ij}$ a {\em restriction} for the color $C_j$.

In \cite{g&j} Greenwell and Johnson define $\X k S$ to be the smallest
integer $m$ such that $S$ is $D$-colorable for every $k \by m$ array $D$
of positive real numbers.  We call $m$ the {\em$k^{\rm th}$ upper
chromatic number} of $S$.  
We allow the case $k=0$, and note that $\X 0 S =1$ for any (nonempty) $S$.  
(In the case $k=1$, we call $m$ the 
{\em upper chromatic number} of $S$, and write simply $\Xone S$, as in 
\cite{g&j}.)
If no such integer exists then we say that $\X k S = \infty$.  Note
that $\X k S$ is a generalization of the ordinary chromatic number $\chi
(S)$, for which only the $1\by m$ array consisting of all 1's is considered.

Very little is known about $\X k S$, even in the case $k=1$.  Greenwell
and Johnson \cite{g&j}
obtained Lemma 1 and posed Problem 1, which is still open.

\begin{lem}[Greenwell, Johnson] \label{XoneR}
$\Xone \R =3$.
\end{lem}

\begin{prob}
Is $\Xone {\R^2}$ finite?
\end{prob}

\section{Lower Bounds} \label{lowerbounds}

In \cite{g&j} it is observed that $\X k \Z > k,$ since $S=\{0,1,\ldots,
k\}$ cannot be $D$-colored for the $k\by k$ array $D=
\pmatrix{
1 & 1 & \ldots & 1 \cr
2 & 2 & \ldots & 2 \cr
\vdots & \vdots & \ddots & \vdots \cr
k & k & \ldots & k \cr
}.$  
(Therefore $\X k \R > k$ as well.)  Indeed, \X k {\R^n} must be greater
than the cardinality of any $k$-distance set in $\R^n$.  (A {\em 
$k$-distance set\/} is a set
of points in $\R^n$ whose pairwise distances take on at most $k$ values.)
Let $m(n,k)$ denote the maximum size of a $k$-distance set in $\R^n$.  The
best general lower bound for $m(n,k)$ (for $n$ large) is $m(n,k) \geq {n+1
\choose k}$ (see for example \cite{laci}, Exercise 1.2.21).  This gives 
us the following theorem.
\begin{thm}  \label{lowbound}
$\X k {\R^n} \geq {n+1\choose k}. \quad\natural$
\end{thm}

In many special cases we can use $k$-distance sets to obtain stronger
results.  For example, $\X k {\R^2} > 2k$ for each $k$, since the vertices
of a regular $(2k+1)$-gon form a $k$-distance set in the plane.
In general, though, we suspect the bound of Theorem \ref{lowbound} to be
very poor.  Next, we present our only improvement on the lower bound in
one dimension, which does not directly use $k$-distance sets.


\begin{thm}  $\X 2 \Z \geq 4$.
\end{thm}
{\bf Proof.}  Let $D={1\,1\,1\choose2\,3\,4}$
and suppose that $\Z$ is $D$-colorable.  Note that $\Z$ is not 
$1\,1\choose3\,4$-colorable;
hence a $D$-coloring of $\Z$ must contain a point colored $C_1$.  Assume
0 is colored $C_1$.  Then $\pm 1,\pm 2$ cannot be colored $C_1$.  Also,
2 and $-2$ cannot both be colored $C_3$; assume 2 is colored $C_2$.
Then 1 must be colored $C_3$, and 3 must be colored $C_1$ or $C_3$ 
(see figure).  If
3 is colored $C_1$ then we cannot color 5, but if 3 is colored $C_3$
then we cannot color -1. $\quad\natural\\$
$$\matrix{
-2 & -1 & 0 & 1 & 2 & 3 & 4 & 5 \cr
\bullet & \bullet & \bullet & \bullet & \bullet & \bullet & \bullet & \bullet & \cr
\phantom{C_3} & \phantom{C_3} & C_1 & C_3 & C_2 & ? & \phantom{C_3} & \phantom{C_3} \cr }$$
\\
\\
The preceding proof can be modified to show that $\Z$ is not $D$-colorable
for $$D=\pmatrix{
1 & 1 & 1 \cr
2n & 2n+1 & 2n+2 \cr },$$
where $n \in \N$.  It might be interesting to investigate the question:  
For which $k \by m$ arrays $D$ is $S\ D$-colorable ($m < \X k S$), for $S=
\Z, \R, \R^n$, and so on.

\section{Upper Bounds}

We now restrict our attention to one dimension, as we head toward an 
upper bound for \X k \R.
The following theorem will be essential in establishing this bound.


\vskip 14 pt
\noindent{\bf Lov\' asz Local Lemma \cite{lll}}
{\em Let $A_1,\ldots,A_n$ be (bad) events in a probability space $\Omega$, with
the probabilities $\Pr(A_i) < p$ for each $i$.  Let $G$ be any (fixed) dependency graph for $\{ A_i\}$; that is,
$V(G)=\{A_i\}$ and for each
$A_i \in V(G)$, the event $A_i$ is independent of the set $\{A_j: A_i A_j
\not\in E(G)\}$.  Let $\D$ be the maximum degree of $G$.  
If $4p\D < 1$, then $\Pr(\cap A_i^c) > 0$.} $\quad\natural\\$


For a $k\by m$ array $D$ (of positive real numbers, as always), 
define $\r(D)$ to be the largest number of columns of $D$ 
in which any particular restriction appears.  Thus, for example, 
$$\r\pmatrix{
1 & 1 & 2 & 3 \cr
1 & 1 & 4 & 5 \cr
1 & 6 & 6 & 6 \cr}
=3,$$ since ``6'' appears in 3 columns.


Now, let $S$ be a finite subset of $\R$, and let $D$ be a $k\by m$ array.  Let
$\Omega$ be the set of all colorings $T:S\to\{\colors\}$.  For each $x,y
\in S$, let $A_{xy}$ be the event that ``for some $i$, $x$ and $y$ are
both colored $C_i$ and $|x-y|$ is a restriction for $C_i$.''  Then $\Pr
(A_{xy})$ is at most $p={\r (D) \over m^2}$, since $|x-y|$ is restricted
for at most $\r (D)$ colors.  Also, the event $A_{xy}$ affects only those
events $A_{xz}$ where $|x-z| \in D$ and $A_{wy}$ where $|w-y| \in D$.
Since there are at most $km-\r (D)+1$ distinct restrictions in $D$, at
most $2(km-\r (D)+1)$ numbers $z$ satisfy $|x-z|\in D$ and at most $2(
km-\r (D)+1)$ numbers $w$ satisfy $|w-y|\in D$.  Furthermore, one such
value of $z$ is $z=y$, and one such value of $w$ is $w=x$.  Hence there
is a dependency graph $G$ with
\begin{eqnarray*}
\D & \leq & 4(km-\r (D)+1) -2 \\
& < & 4km,
\end{eqnarray*}
and
\begin{eqnarray*}
4p\D & < & 4\cdot {\r (D) \over m^2} \cdot 4km \\
& = & {16k\r (D) \over m}. 
\end{eqnarray*}
Thus if $m\geq 16k\r (D)$, then $4p\D <1$, and by the Lov\' asz Local
Lemma, $\cap A_{xy}^c$ has
positive probability.  That is, there is some coloring of $S$ obeying
all restrictions, so $S$ is $D$-colorable.  We have just proven the 
following lemma.
\begin{lem}  \label{finite}
Let $S$ be a finite subset of $\R$, and let $D$ be a $k\by m$
array.  If $16k\r (D) \leq m$, then $S$ is $D$-colorable. 
\end{lem}


The next lemma is an extension of Lemma \ref{finite} to the
case where $S$ is countable.
\begin{lem}  \label{countable}
Let $S=\{x_1,x_2,\ldots\}\subseteq\R$, and let $D$ be a 
$k\by m$ array.  If $16k\r (D) \leq m$, then $S$ is $D$-colorable.
\end{lem}
{\bf Proof.}  For each $n$, let $S_n=\{x_1,\ldots,x_n\}$.  Assume $16k\r (D)
\leq m$.  By Lemma \ref{finite}, $S_n$ is $D$-colorable for each $n$.
Let $T_n$ be a $D$-coloring of $S_n$, and consider the
sequence $T_1(x_1)$, $T_2(x_1)$, $\cdots$ of colors.  Since we have
only $m$ colors, this sequence must contain an infinite monochromatic
subsequence $T_{i_1}(x_1)$, $T_{i_2}(x_1)$, $\cdots$.  Color $x_1$
with the color in this subsequence.  Now consider the sequence $T_{i
_1}(x_2)$, $T_{i_2}(x_2)$, $\cdots$.  This too must contain an infinite
monochromatic subsequence; color $x_2$ with the color of this
subsequence.  Continuing in this manner, we obtain a $D$-coloring of
$S$, and the proof is complete.
\qed


Finally, we extend to the general case, $S=\R$.
\begin{prop}  \label{smallrho}
Let $D$ be a $k \by m$ array.  If $16k\r (D) \leq m$, then $\R$ is
$D$-colorable.
\end{prop}
{\bf Proof.}  Let $G$ be the (additive) subgroup of $\R$ generated by
the elements of $D$.  Since $G$, and therefore each coset of $G$, is
countable, we can $D$-color any coset of $G$ by Lemma \ref{countable}. 
Thus we can color $\R$ by $D$-coloring every coset of $G$.  We claim that
this gives a $D$-coloring of $\R$.  For, let $x,y \in \R$.  If $x-y \in
G$, then no restriction is violated, since we have $D$-colored the coset
$x+G$ ($=y+G$).  If $x-y \not\in G$, then in particular $|x-y| \not\in D$,
so again no restriction is violated.
\qed
%
%
Proposition \ref{smallrho} tells us that we can $D$-color $\R$ as long as 
$D$ has sufficiently few repetitions.  We now treat the case in which 
$\r (D)$ is large, starting with $k=2$.  
Let $D$ be a $2 \by m$ array, and assume that the restriction $r$ appears
in columns $1,2,\ldots,\r(D)$ of $D$.  Let the other restrictions in 
these columns be denoted by $a_1,a_2,\ldots,a_{\r(D)}$, so that
$$D=\pmatrix{
r & r & \cdots & r & \cdots \cr
a_1 & a_2 & \cdots & a_{\r(D)} & \cdots \cr
}.$$
\vskip 28 pt
\begin{lem} \label{bigrho2}
Let $D$ be a $2\by m$ array.  If $\r (D) \geq 6$ then $\R$ is $D$-colorable.
\end{lem}
{\bf Proof.}  Divide $\R$ into half-open intervals $[nr,(n+1)r)$, for each $n
\in \Z$.  We define such an interval to be {\em even} if $n$ is even and
{\em odd} if $n$ is odd.  Let $U$ be the union of the even intervals, and
let $V$ be the union of the odd intervals.  
Since $\Xone \R =3$ (Lemma \ref{XoneR}), we can color $U$ with the colors 
$C_1, C_2, C_3$ subject only to the restrictions $a_1, a_2, a_3$, 
respectively.  Similarly, we can color $V$ with the colors $C_4, C_5, C_6$,
subject to the restrictions $a_4, a_5, a_6$.  Thus we can combine these
colorings to get a coloring of $U \cup V = \R$.  This coloring is in fact
a $D$-coloring of $\R$, since clearly no two points colored $C_i$
are a distance $r$ apart.
\qed
This idea extends quite easily to all $k$, as follows.  
%
%
\begin{lem} \label{bigrho}
Let $D=(d_{ij})$
be a $k\by m$ array with $\r (D)=\r$, and assume that $d_{11}=d_{12}=
\cdots=d_{1\r}=r$.  Let $A$ be the $(k-1)\by\r$ array obtained
from the first $\r$ columns of $D$ by deleting the top row.
Assume $\X {k-1} \R$ is finite and $\r\geq 2\X {k-1} \R$.  Then $\R$ is 
$D$-colorable.
\end{lem}
{\bf Proof.}  Let $U$ and $V$ be as they were in Lemma \ref{bigrho2}.
Because of our assumption regarding $\r$, we may define
$A'$ to be the array consisting of the first $\X {k-1} \R$ columns
of $A$, and $A''$ to be the array consisting of the next $\X {k-1} \R$ 
columns of $A$.  Since $A$
has just $k-1$ rows, we can $A'$-color $U$ and $A''$-color $V$.  Combining
the colorings gives an $A$-coloring of $\R$, which is also a $D$-coloring
since no two points colored $C_i$ are a distance $r$ apart.
\qed
%
%
\begin{thm} \label{XkR}
For each $k$, $\X k \R$ is finite.
\end{thm}
{\bf Proof.}  We proceed by induction.  The case $k=0$ is trivial,
and $k=1$ is Lemma \ref{XoneR}.  Assume $\X {k-1} \R$ is finite, and let
$$m=32k\X {k-1} \R -16k.$$
We will show that $\X k \R \leq m$.

Let $D$ be a $k \by m$ array.  If $\r (D) < 2\X {k-1} \R$ then $\r (D)
\leq 2\X {k-1} \R -1$, so $\R$ is $D$-colorable by Proposition 
\ref{smallrho}.  If $\r (D) \geq 2\X {k-1} \R$ then $\R$ is $D$-colorable 
by Lemma \ref{bigrho}.  A third case is impossible. 
\qed


In particular, we have $\X k \R \leq 32^k \cdot k!$.
The bound of Theorem \ref{XkR} is the best upper bound we have 
obtained, other than by using a refined version of the Lov\' asz
Local Lemma. 
Note that the existence of this bound does not depend on Lemma
\ref{XoneR}:  from Proposition \ref{smallrho} it follows that
$\Xone \R \leq 16$.  (For any
$1\by m$ array $D$, if $\r(D)>1$ then $\R$ is $D$-colorable.)

\section{Conclusions and Further Problems}

We have shown that $\X k \R$ is finite for all $k$.  However, our
methods completely fail to generalize to higher dimensions.  Thus
the most general open question is still:
\begin{prob}
Is $\X k {\R^n}$ finite for all $k$ and $n$?
\end{prob}
It would also be nice to know how $\X k {\R^n}$ grows, as a function
of $k$ and $n$.  
Even the behaviour of $\X k \Z$ is mysterious, however.  In each of the
dozens of examples we have examined, we observe that whenever $\Z$
is $D$-colorable, we can find a periodic $D$-coloring of $\Z$;
we wonder if this is always the case.  Next, it is easy to see
that for each $k$, $\X k \Z = \X k \Q$.  Perhaps this motivates
\begin{prob}
Is $\X k \R = \X k \Z$?
\end{prob}
Finally, as suggested by Section \ref{lowerbounds}, we pose the perhaps
easier question:
\begin{prob}
Is $\X 2 \Z =4$?
\end{prob}


\bigskip
\medskip
\noindent
{\Large\bf Acknowledgements}
\\
\\
The author thanks Rob Hochberg for
introducing him to the Lov\' asz Local Lemma, and Joe Gallian
and the University of Minnesota, Duluth for making his research
possible.  He also thanks The Geometry Center at the University 
of Minnesota for the use of their computers, and several reviewers 
for their helpful comments and suggestions, especially Dean 
Hickerson.  Work supported by NSF Grant \#DMS9000742 and 
NSA Grant \#MDA904-H-0036.

\end{document}